\documentclass[journal]{IEEEtran}
\usepackage{graphicx}
\usepackage{amscd}
\usepackage[cmex10]{amsmath}
\usepackage{amsfonts}
\usepackage{amssymb}
\usepackage{latexsym}
\usepackage{mathrsfs}
\usepackage{fixltx2e}
\usepackage{epstopdf}

\newtheorem{theorem}{Theorem}

\newtheorem{corollary}{Corollary}[theorem]

\renewcommand{\IEEEQED}{\IEEEQEDopen}
\begin{document}
\title{Monotonicity properties of some Dini functions}

\author{\'Arp\'ad~Baricz, Tibor~K.~Pog\'any and R\'obert Sz\'asz
\thanks{{\bf \'Arp\'ad Baricz} is with Department of Economics, Babe\c{s}-Bolyai University, 400591 Cluj--Napoca, Romania and Institute of Applied Mathematics, John von Neumann Faculty of Informatics, \'Obuda University, 1034 Budapest, Hungary; {\bf Tibor K. Pog\'any} is with Faculty of Maritime Studies, University of Rijeka, 51000 Rijeka, Croatia and Institute of Applied Mathematics, John von Neumann Faculty of Informatics, \'Obuda University, 1034 Budapest, Hungary; {\bf R\'obert Sz\'asz} is with Department of Mathematics and Informatics, Sapientia Hungarian University of Transylvania, 540485 T\^argu Mure\c{s}, Romania.}%
\thanks{{\bf e--mail:}~{bariczocsi@yahoo.com} (\'A. Baricz),~{tkpogany@gmail.com} (T.K. Pog\'any),~{rszasz@ms.sapientia.ro} (R. Sz\'asz)}}

\markboth{IEEE 9th International Symposium on Applied Computational
Intelligence and Informatics}%
{}

\maketitle

\begin{abstract}
In this note our aim is to deduce some new monotonicity properties for a special combination of Bessel functions of the first kind by using a recently developed Mittag-Leffler expansion for the derivative of a normalized Bessel function of the first kind. These monotonicity properties are used to obtain some new inequalities for Bessel functions of the first kind.
\end{abstract}


\section{Introduction}

\IEEEPARstart{B}{essel} functions of the first kind play an important role in various branches of applied mathematics and engineering sciences. Their properties have been investigated by many scientists and there is a very extensive literature dealing with Bessel functions. One of the most important properties of Bessel functions of the first kind is that they can be represented by an Weierstrassian infinite product. This is a very useful tool in proving many inequalities and important properties of Bessel functions. Recently, there has been a vivid interest on geometric properties of Bessel functions, like univalence, starlikeness, convexity, close-to-convexity in the open unit disk. For more details we refer to the paper \cite{szasz} and to the references therein. In \cite{szasz}, in order to determine the radius of convexity of a normalized Bessel function, Baricz and Sz\'asz proved a Mittag-Leffler expansion for a special combination of Bessel functions, called sometimes as Dini function. In this note our aim is to deduce an infinite product representation for this Dini function and to apply this result in order to generate some absolute monotonicity properties of some quotients of Dini functions. This paper is organized as follows: in the next section we present the above mentioned infinite product representation of Dini functions. Section III contains the monotonicity properties together with inequalities, while Section IV contains the concluding remarks.

\section{Infinite product representations of some Dini functions}

The following result is the key tool in the proof of our main results of Section III. The proof of the next infinite product representation is based on a Mittag-Leffler expansion, recently developed by Baricz and Sz\'asz \cite[Lemma 4]{szasz}, however, we present an alternative proof by using the Hadamard theorem concerning the growth order of entire functions. Throughout in the sequel $J_{\nu}$ denotes the Bessel function of the first kind of order $\nu.$

\begin{theorem}
{\em Let $\nu>-1$ and consider the Dini function $d_{\nu}:\mathbb{C}\to\mathbb{C},$ defined by $$d_{\nu}(z)=(1-\nu)J_{\nu}(z)+zJ_{\nu}'(z).$$ If $\alpha_{\nu,n}$ denotes the $n$th positive zero of the Dini function $d_{\nu},$ then the following Weierstrassian factorization is valid
\begin{equation}\label{product}
d_{\nu}(z)=\frac{z^{\nu}}{2^{\nu}\Gamma(\nu+1)}\prod_{n\geq1}\left(1-\frac{z^2}{\alpha_{\nu,n}^2}\right),
\end{equation}
where the infinite product is uniformly convergent on each compact subset of the complex plane.}
\end{theorem}

\noindent {\it Proof:} Let $z$ be a complex number. Consider the function $g_{\nu},$ defined by $$g_{\nu}(z)=2^{\nu}\Gamma(\nu+1)z^{1-\nu}J_{\nu}(z).$$
We know that for $\nu>-1$ we have \cite[Lemma 4]{szasz}
\begin{equation}\label{log}\frac{g_{\nu}''(z)}{g_{\nu}'(z)}=-\sum_{n\geq1}\frac{2z}{\alpha_{\nu,n}^2-z^2}.\end{equation}
Integrating both sides of the above relation we get
$$\log g_{\nu}'(z)=\sum_{n\geq 1}\log\left(1-\frac{z^2}{\alpha_{\nu,n}^2}\right)+c,$$
which implies
$$g_{\nu}'(z)=e^{c}\prod_{n\geq1}\left(1-\frac{z^2}{\alpha_{\nu,n}^2}\right),$$
where $c$ is a constant. On the other hand, we have
$$g_{\nu}'(z)=2^{\nu}\Gamma(\nu+1)z^{-\nu}d_{\nu}(z),$$
which by using the recurrence relation
\begin{equation}\label{rec}zJ_{\nu}'(z)=-zJ_{\nu+1}(z)+\nu J_{\nu}(z),\end{equation}
can be rewritten as
\begin{equation}\label{rec2}g_{\nu}'(z)=2^{\nu}\Gamma(\nu+1)z^{-\nu}\left(J_{\nu}(z)-zJ_{\nu+1}(z)\right).\end{equation}
Now, by using the infinite sum representation of the Bessel function of the first kind
\begin{equation}\label{sum}J_{\nu}(z)=\sum_{n\geq 0}\frac{(-1)^nz^{2n+\nu}}{2^{2n+\nu}n!\Gamma(n+\nu+1)}\end{equation}
we obtain $g_{\nu}'(0)=1,$ which in turn implies that $c=0.$ This completes the proof of the infinite product representation \eqref{product}. Finally, we mention that the uniform convergence of the infinite product in \eqref{product} is a consequence of the uniform convergence of the infinite sum in \eqref{log}, which follows by using the proof of \eqref{log}, see \cite[Lemma 4]{szasz}.

Alternatively, the formula \eqref{product} can be proved as follows. From \eqref{rec2} and \eqref{sum} we obtain that
$$g_{\nu}'(z)=1+\sum_{n\geq0}\frac{2n+3}{2n+2}\cdot\frac{(-1)^n \Gamma(\nu+1)z^{2n+2}}{2^{2n+1}n!\Gamma(n+\nu+2)}.$$
Taking into consideration the well-known limits
$$
\lim_{n\to \infty} \frac{\log \Gamma(n+c)}{n\,\log n}=1,\quad
\lim_{n\to \infty} \frac{[\Gamma(n+c)]^{1/n}}{n}=\frac{1}{e},
$$
where $c$ is a positive constant, and \cite[p. 6, Theorems 2 and 3]{lev}, we infer that the entire
function $g_{\nu}'$ is of growth order $\rho=\frac{1}{2}$ and exponential type $\sigma=1.$ Namely,
for $\nu>-1$ we have that as $n\to\infty$
$$\frac{n\log n}{\log\frac{2^{2n+1}\Gamma\left(n+1\right)}{\Gamma(\nu+1)}+\log\Gamma\left(n+\nu+2\right)+\log\left(\frac{2n+2}{2n+3}\right)}
\to\frac{1}{2}$$
and
$$\frac{n}{\rho e}\sqrt[2n]{\frac{2n+3}{2n+2}\cdot\frac{\Gamma(\nu+1)}{2^{2n+1}\Gamma(n+1)\Gamma(n+\nu+2)}}\to 1.$$
Now, recall that \cite[p. 597]{watson} in case $\alpha+\nu>0$
and $\nu>-1$ the Dini function $z\mapsto zJ_{\nu}'(z)+\alpha J_{\nu}(z)$ has only real zeros. With this the rest of the proof of \eqref{product} follows by applying Hadamard's Theorem \cite[p. 26]{lev}. \medskip \hfill $\Box$

\section{Monotonicity properties of some combinations of Bessel functions of the first kind}

Our first main result of this paper is the following theorem, which presents the absolute monotonicity of three functions involving the Dini function $d_{\nu}$. The proofs borrow some ideas from the paper of Ismail and Muldoon \cite{isma}. Note that a function $f:I\to\mathbb{R}$ is called absolutely monotonic on the interval $I$ if for all $x\in I$ and $n\in\{0,1,\dots\}$ we have $f^{(n)}(x)\geq0.$

\begin{theorem}
{\em Let $\mu\geq \nu>-1$ and consider the Dini function $d_{\nu}:\mathbb{R}\to\mathbb{R},$ defined by $$d_{\nu}(x)=(1-\nu)J_{\nu}(x)+xJ_{\nu}'(x).$$ Then the functions $f_{\mu,\nu},g_{\mu,\nu},q_{\nu}:[0,\alpha_{\nu,1}^2)\to(0,\infty),$ defined by
\begin{align*}&f_{\mu,\nu}(x)=\left(\log\left(x^{\frac{\nu-\mu}{2}}e^{\frac{3x}{4}\left(\frac{1}{\mu+1}-\frac{1}{\nu+1}\right)}
\frac{d_{\mu}(\sqrt{x})}{d_{\nu}(\sqrt{x})}\right)\right)',\\
&g_{\mu,\nu}(x)=x^{\frac{\nu-\mu}{2}}e^{\frac{3x}{4}\left(\frac{1}{\mu+1}-\frac{1}{\nu+1}\right)}
\frac{d_{\mu}(\sqrt{x})}{d_{\nu}(\sqrt{x})},\\
&q_{\nu}(x)=\frac{x^{\frac{\nu}{2}}e^{-\frac{3x}{4(\nu+1)}}}{d_{\nu}(\sqrt{x})},\end{align*}
are absolutely monotonic.}
\end{theorem}

\noindent{\it Proof:} By using \eqref{product} we get
$$\left(\log\left( x^{-\frac{\nu}{2}}d_{\nu}(\sqrt{x})\right)\right)'=\sum_{n\geq1}\frac{1}{x-\alpha_{\nu,n}^2}$$
and hence
$$f_{\mu,\nu}(x)=\sum_{n\geq1}\left(\frac{1}{\alpha_{\nu,n}^2-x}-\frac{1}{\alpha_{\mu,n}^2-x}\right)+\frac{3(\nu-\mu)}{4(\nu+1)(\mu+1)}.$$
Recall that \cite[p. 196]{landau} if $\gamma_{\nu,n}$ is the $n$th positive root of the equation $$\gamma{J}_\nu(z)+zJ_\nu'(z)=0$$ and $\nu+\gamma\geq0,$ then the function $\nu\mapsto \gamma_{\nu,n}$ is strictly increasing on $(-1,\infty)$ for $n\in\{1,2,\dots\}$ fixed. Consequently, we have that $\nu\mapsto \alpha_{\nu,n}$ is strictly increasing on $(-1,\infty)$ for $n\in\{1,2,\dots\}$ fixed. Thus, for all $n,m\in\{1,2,\dots\},$ $\mu\geq\nu>-1$ and $x\leq \alpha_{\nu,1}^2$ we have $$(\alpha_{\nu,n}^2-x)^{m+1}\leq(\alpha_{\mu,n}^2-x)^{m+1},$$ and this in turn implies that
$$f_{\mu,\nu}^{(m)}(x)=\sum_{n\geq 1}\left(\frac{m!}{(\alpha_{\nu,n}^2-x)^{m+1}}-\frac{m!}{(\alpha_{\mu,n}^2-x)^{m+1}}\right)\geq0$$
for all $\mu\geq\nu>-1,$ $m\in\{1,2,\dots\}$ and $x\leq \alpha_{\nu,1}^2.$

Now, dividing each side of \cite[Lemma 4]{szasz}
$$\frac{zJ_{\nu+2}(z)-3J_{\nu+1}(z)}{J_{\nu}(z)-zJ_{\nu+1}(z)}
=-\sum_{n\geq1}\frac{2z}{\alpha_{\nu,n}^2-z^2}$$
with $z$ and tending with $z$ to $0,$ in view of \eqref{sum} we obtain that
$$\sum_{n\geq 1}\frac{1}{\alpha_{\nu,n}^2}=\frac{3}{4(\nu+1)}.$$
This implies that $f_{\mu,\nu}(0)=0,$ and by using the fact that $f_{\mu,\nu}$ is increasing we get that $f_{\mu,\nu}(x)\geq0$ for all $\mu\geq\nu>-1$ and $x\in[0,\alpha_{\nu,1}^2).$ This completes the proof of the absolute monotonicity of $f_{\mu,\nu}.$

Since the exponential of a function having an absolutely monotonic
derivative is absolutely monotonic, it follows that $g_{\mu,\nu}$ is also absolutely monotonic on $[0,\alpha_{\nu,1}^2).$

Finally, we note that by using \eqref{product} it follows that as $\mu\to \infty$ we have
$$2^{\mu}\Gamma(\mu+1)x^{-\mu}d_{\mu}(x)\to 1,$$ which in turn implies that the absolute monotonicity of $q_{\nu}$ follows from the absolute monotonicity of $g_{\mu,\nu}$ by tending with $\mu$ to infinity. \medskip \hfill $\Box$

We note that the above theorem can be used to prove inequalities for the Dini function $d_{\nu}.$ For example, the absolute monotonicity of $q_{\nu}$ in particular implies the following result.

\begin{corollary}
{\em If $\nu>-1$ and $x\in[0,\alpha_{\nu,1}),$ then
\begin{equation}\label{bound}d_{\nu}(x)\leq\frac{x^{\nu}e^{-\frac{3x^2}{4(\nu+1)}}}{2^{\nu}\Gamma(\nu+1)}.\end{equation}}
\end{corollary}

\noindent {\it Proof:} By using again \eqref{product} we obtain $q_{\nu}(0)=2^{\nu}\Gamma(\nu+1)$ and since $q_{\nu}$ is increasing on $[0,\alpha_{\nu,1}^2)$ it follows that $q_{\nu}(x)\geq q_{\nu}(0),$ that is,
$$2^{\nu}\Gamma(\nu+1)d_{\nu}(\sqrt{x})\leq x^{\frac{\nu}{2}}e^{-\frac{3x}{4(\nu+1)}}.$$
Changing $x$ to $x^2$ the proof is complete. \medskip \hfill $\Box$

\begin{figure}[htb] \label{plot}
      \centering
        \scalebox{0.65}
        {\includegraphics{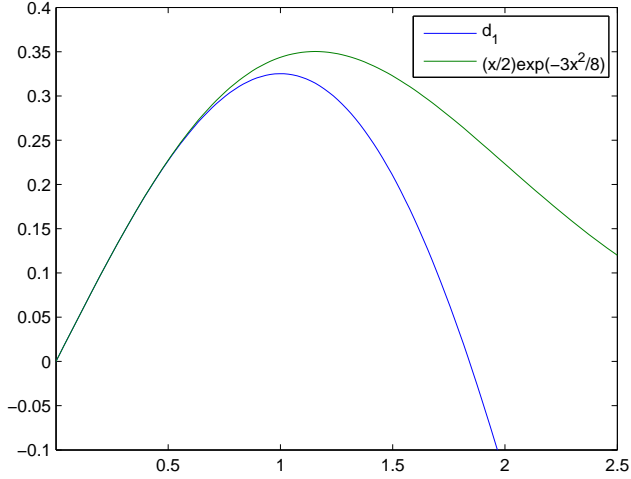}}
        \caption{Graph of the functions $d_1$ and $x\mapsto \frac{x}{2}e^{-\frac{3x^2}{8}}$.}
\end{figure}

It is worth to mention here that by using the similar inequality for the Bessel function of the
first kind $J_{\nu}$ we can obtain a similar inequality for $d_{\nu}$ as above. Namely, we know that \cite[Corollary 7]{isma} for $\nu>-1$ the function $$x\mapsto \frac{x^{\frac{\nu}{2}}e^{-\frac{x}{4(\nu+1)}}}{J_{\nu}(\sqrt{x})}$$
is absolutely monotonic on $[0,j_{\nu,1}^2)$ and hence we have
\begin{equation}\label{bound1}
J_{\nu}(x)\leq \frac{x^{\nu}e^{-\frac{x^2}{4(\nu+1)}}}{2^{\nu}\Gamma(\nu+1)}
\end{equation}
for $\nu>-1$ and $x\in[0,j_{\nu,1}).$ Here $j_{\nu,1}$ denotes the first positive zero of the Bessel function $J_{\nu}.$ By using \eqref{bound1} and the fact that $J_{\nu+1}(x)>0$ for $x\in(0,j_{\nu+1,1})$ we get for $\nu>-1$ and $x\in[0,j_{\nu,1})$
$$d_{\nu}(x)=J_{\nu}(x)-xJ_{\nu+1}(x)\leq \frac{x^{\nu}e^{-\frac{x^2}{4(\nu+1)}}}{2^{\nu}\Gamma(\nu+1)},$$
however, this upper bound is not better than the upper bound in \eqref{bound}.

Note also that Fig. 1 suggest that the lower bound in \eqref{bound} is quite tight near the origin. To prove this observe that from \eqref{product}
we have that $$g_{\nu}'(x)=2^{\nu}\Gamma(\nu+1)x^{-\nu}d_{\nu}(x)$$ takes the value $1$ at the origin, and we have also that $e^{-\frac{3x^2}{4(\nu+1)}}$ becomes $1$ as $x=0.$ These show that indeed the inequality \eqref{bound} is sharp as $x=0.$

Now, we are going to show another two applications of the infinite product representation \eqref{product}. In what follows we use the notations $\alpha_{\nu,0}=0$ and $$\Delta_{\nu}=\bigcup_{k\geq0}(\alpha_{\nu,2k};\alpha_{\nu,2k+1}).$$
The first result is about the log-concavity of the function $d_{\nu},$ while the second deals with the log-concavity of $g_{\nu}',$ which in particular yields a new van der Corput type inequality for Bessel functions. For more details on van der Corput inequalities for Bessel and modified Bessel functions of the first kind we refer to the paper \cite{blp} and to the references therein.

\begin{theorem}
{\em If $\nu>-1,$ then the function $f_{\nu}:\Delta_{\nu}\to\mathbb{R},$ defined by $$f_{\nu}(x)=\frac{d_{\nu}'(x)}{d_{\nu}(x)}-\frac{\nu}{x},$$
is decreasing and consequently the function $d_{\nu}$ is log-concave on $\Delta_{\nu}$ for $\nu\geq0.$ Consequently, we have
for all $\nu>-1$ and $x\in\Delta_{\nu}$ the inequality
\begin{equation}\label{logar}\frac{xd_{\nu}'(x)}{d_{\nu}(x)}<\nu.\end{equation}}
\end{theorem}

\noindent {\it Proof:} First note that by using \eqref{product} we have that $d_{\nu}(x)>0$ for $x\in\Delta_{\nu},$ and thus the log-concavity of $d_{\nu}$
on $\Delta_{\nu}$ makes sense. By using logarithmic differentiation from \eqref{product} we get
$$f_{\nu}(x)=-\sum_{n\geq1}\frac{2x}{\alpha_{\nu,n}^2-x^2},$$
which implies that
$$f_{\nu}'(x)=-2\sum_{n\geq 1}\frac{\alpha_{\nu,n}^2+x^2}{(\alpha_{\nu,n}^2-x^2)^2}<0$$
for all $x\in\Delta_{\nu}$ and $\nu>-1.$ Thus, we have
$$\left(\frac{d_{\nu}'(x)}{d_{\nu}(x)}\right)'<-\frac{\nu}{x^2}\leq0$$
for all $\nu\geq 0$ and $x\in\Delta_{\nu},$ and consequently the function $d_{\nu}$ is indeed log-concave on $\Delta_{\nu}$ for $\nu\geq0.$
On the other hand, by using the above representation of $f_{\nu}$ we obtain $f_{\nu}(0)=0,$ and thus we get $f_{\nu}(x)<0,$ which is equivalent to
\eqref{logar}. \medskip \hfill $\Box$

\begin{theorem}
{\em The function $g_{\nu}'$ is log-concave on $\Delta_{\nu}$ for all $\nu>-1$ and consequently for all $\nu>-1$ and $a,b\in\Delta_{\nu}$ the following van der Corput type inequality is valid
\begin{equation}\label{corput}\left|g_{\nu}(a)-g_{\nu}(b)\right|\geq |a-b|\sqrt{g_{\nu}'(a)g_{\nu}'(b)}\end{equation}
or equivalently
$$\left|a^{1-\nu}J_{\nu}(a)-b^{1-\nu}J_{\nu}(b)\right|\geq |a-b|\sqrt{(ab)^{-\nu}d_{\nu}(a)d_{\nu}(b)}.$$}
\end{theorem}

\noindent {\it Proof:} By using \eqref{product} we get
$$\left(\frac{g_{\nu}''(x)}{g_{\nu}'(x)}\right)'=\left(-\sum_{n\geq 1}\frac{2x}{\alpha_{\nu,n}^2-x^2}\right)'=-\sum_{n\geq 1}\frac{2(\alpha_{\nu,n}^2+x^2)}{(\alpha_{\nu,n}^2-x^2)^2}$$
and thus the function $g_{\nu}'$ is log-concave on $\Delta_{\nu}$ for all $\nu>-1.$

Now, we shall use the following result. Let $f:[a,b]\to\mathbb{R}$ be such that $f':[a,b]\to(0,\infty)$ is log-concave. Then the inequality \cite[p. 242]{corput}
$$\log\left(\frac{1}{b-a}\int_a^bf'(x)dx\right)\geq \frac{1}{b-a}\int_a^b\log f'(x)dx$$
$$\geq\frac{\log f'(a)+\log f'(b)}{2}$$
implies that
\begin{equation}\label{fcorput}\log \frac{f(b)-f(a)}{b-a}\geq \log\sqrt{f'(a)f'(b)}\end{equation}
is valid. Applying \eqref{fcorput} for the function $f=g_{\nu}$ we get \eqref{corput}. \medskip \hfill $\Box$

By using the relations
$$J_{-\frac{1}{2}}(x)=\sqrt{\frac{2}{\pi x}}\cos x, \ \ J_{\frac{1}{2}}(x)=\sqrt{\frac{2}{\pi x}}\frac{\sin x}{x},$$
the above result in particular reduces to the following.

\begin{corollary}
{\em Let $\alpha_{-\frac{1}{2},n}$ be the $n$th positive zero of $d_{-\frac{1}{2}},$ that is, the $n$th positive root of the equation
$\cos x=\sin x.$ Then the following van der Corput type inequality is valid
$$\left|a\cos a-b\cos b\right|\geq |a-b|\sqrt{(\cos a-\sin a)(\cos b-\sin b)},$$
where $$a,b\in\Delta_{-\frac{1}{2}}=\left(0,\frac{\pi}{4}\right)\bigcup\left(\frac{5\pi}{4},\frac{9\pi}{4}\right)\bigcup{\dots}.$$}
\end{corollary}

Now, for $m\in\{1,2,\dots\}$ let us use the following notation
$$\eta_{2m}(\nu)=\sum_{n\geq 1}\frac{1}{\alpha_{\nu,n}^{2m}}.$$
Another interesting applications of the infinite product in \eqref{product} are the following results.

\begin{theorem}
{\em If $\nu>-1$ and $z\in\mathbb{C}$ such that $|z|<\alpha_{\nu,1},$ then we have the following power series representation
\begin{equation}\label{quo}\frac{zd_{\nu}'(z)}{d_{\nu}(z)}=\nu-2\sum_{m\geq 1}\eta_{2m}(\nu)z^{2m}.\end{equation}
In particular, the function
$$x\mapsto -\frac{xd_{\nu}'(x)}{d_{\nu}(x)}=-\nu+2\sum_{m\geq 1}\eta_{2m}(\nu)x^{2m}$$
is absolutely monotonic on $(0,\alpha_{\nu,1})$ for all $\nu>-1.$}
\end{theorem}

\noindent {\it Proof:} By using again \eqref{product} we have
\begin{align*}\frac{zd_{\nu}'(z)}{d_{\nu}(z)}&=\nu-2\sum_{n\geq 1}\frac{z^2}{\alpha_{\nu,n}^2-z^2}\\&
=\nu-2\sum_{n\geq 1}\sum_{m\geq1}\left(\frac{z^2}{\alpha_{\nu,n}^2}\right)^m\\
&=\nu-2\sum_{m\geq 1}\left(\sum_{n\geq1}\frac{1}{\alpha_{\nu,n}^{2m}}\right)z^{2m},\end{align*}
where $|z|<\alpha_{\nu,1}$ and $\nu>-1.$ \medskip \hfill $\Box$

\begin{theorem}
{\em The functions
$$\nu\mapsto g_{\nu}'(x)=2^{\nu}\Gamma(\nu+1)x^{-\nu}d_{\nu}(x),\ \ \nu\mapsto \frac{xd_{\nu}'(x)}{d_{\nu}(x)}$$
are increasing on $(-1,\infty)$ for all $x\in(0,\alpha_{\nu,1}).$ Equivalently, the following inequalities
are valid for all $x\in(0,\alpha_{\nu,1})$ and $\mu\geq\nu>-1$
\begin{equation}\label{vege}x^{\nu-\mu}\frac{d_{\mu}(x)}{d_{\nu}(x)}\geq 2^{\nu-\mu}\frac{\Gamma(\nu+1)}{\Gamma(\mu+1)}, \ \ \frac{d_{\mu}'(x)}{d_{\mu}(x)}\geq\frac{d_{\nu}'(x)}{d_{\nu}(x)}.\end{equation}
Moreover, the function $x\mapsto d_{\mu}(x)/d_{\nu}(x)$ is increasing on $(0,\alpha_{\nu,1})$ for all $\mu\geq\nu>-1.$}
\end{theorem}

\noindent {\it Proof:} Appealing again on \eqref{product}
we get $$\frac{\partial}{\partial\nu}\left(\log\left(2^{\nu}\Gamma(\nu+1)x^{-\nu}d_{\nu}(x)\right)\right)=
\sum_{n\geq1}\frac{2x^2\frac{\partial}{\partial\nu}\alpha_{\nu,n}}{\alpha_{\nu,n}(\alpha_{\nu,n}^2-x^2)}$$
and since $\nu\mapsto\alpha_{\nu,n}$ is increasing on $(-1,\infty)$ for each $n\in\{1,2,\dots\}$ fixed, it follows that the function
$$\nu\mapsto 2^{\nu}\Gamma(\nu+1)x^{-\nu}d_{\nu}(x)$$ is increasing on $(-1,\infty)$ for $x\in(0,\alpha_{\nu,1}).$

Now, since $\nu\mapsto\alpha_{\nu,n}$ is increasing on $(-1,\infty)$ for each $n\in\{1,2,\dots\}$ fixed, it follows that $\nu\mapsto\eta_{2m}(\nu)$ is decreasing on $(-1,\infty)$ for each $m\in\{1,2,\dots\},$ and by using \eqref{quo}, this implies that $\nu\mapsto \frac{xd_{\nu}'(x)}{d_{\nu}(x)}$
is increasing on $(-1,\infty)$ for all $x\in(0,\alpha_{\nu,1}).$

The inequalities in \eqref{vege} are equivalent with the monotonicity properties, while the last sentence of the theorem is equivalent to the last inequality. With this the proof is complete. \medskip \hfill $\Box$

\section{Conclusion}

In this note we have shown that by using the infinite product representation of a Dini function how we can generate many results on 
the monotonicity properties of quotients of Dini functions. We presented also some inequalities concerning Bessel functions of the first kind. 
As far as we know our results are new and we believe that it would be of interest to study further the properties of the Dini function $d_{\nu}$ in order to deduce some other inequalities for this function.

\end{document}